\documentclass[10pt,twoside]{article}
\usepackage{graphicx}
\usepackage{amsmath,amssymb,amsfonts}
\usepackage{Latex-document}

\newcounter{hack}
\newtheorem{thr}{Theorem}[hack]
\newtheorem{propo}[thr]{Proposition}
\newtheorem{lemm}[thr]{Lemma}
\newtheorem{coro}[thr]{Corollary}
\newtheorem{Defin}[thr]{Definition}
\newtheorem{remark}[thr]{Remark}
\newtheorem{remarks}[thr]{Remarks}
\newtheorem{exam}[thr]{Example}
\newtheorem{exams}[thr]{Examples}

\newenvironment{proof}{\noindent{\bf Proof.}}%
{\hfill$\square$}%
\newenvironment{thm}[1][]{\begin{thr}{$\kern-5pt${\bf #1}{\bf.} }}%
{\end{thr}}
{\end{propo}}
{\end{lemm}}
{\end{coro}}

\newenvironment{Def}{\begin{Defin}$\kern-5pt${\bf.}\rm}%
{\end{Defin}}
\newenvironment{rem}{\begin{remark}$\kern-5pt${\bf.}\rm}%
{\end{remark}}
\newenvironment{rems}{\begin{remarks}$\kern-5pt${\bf.}\rm}%
{\end{remarks}}
{\end{exam}}
\newenvironment{exs}{\begin{exams}$\kern-5pt${\bf.}\rm}%
{\end{exams}}

\newcommand{\CH}{{\rm CH}}

\newcommand{\Pic}{{\rm Pic}}

\newcommand{\Spec}{{\rm Spec \,}}

\newcommand{\sC}{{\mathcal C}}

\newcommand{\sE}{{\mathcal E}}

\newcommand{\sO}{{\mathcal O}}

\newcommand{\sR}{{\mathcal R}}

\newcommand{\sZ}{{\mathcal Z}}
\newcommand{\A}{{\mathbb A}}

\newcommand{\C}{{\mathbb C}}

\renewcommand{\L}{{\mathbb L}}

\newcommand{\N}{{\mathbb N}}
\renewcommand{\P}{{\mathbb P}}
\newcommand{\Q}{{\mathbb Q}}

\newcommand{\U}{{\mathbb U}}

\newcommand{\Z}{{\mathbb Z}}

\newcommand{\id}{{\operatorname{id}}}

\newcommand{\Sch}{{\operatorname{\mathbf{Sch}}}}

\newcommand{\holim}{\mathop{{\rm holim}}}

\newcommand{\Sm}{{\mathbf{Sm}}}

\newcommand{\cn}{{\tilde{c}_1}}

\newcommand{\et}{{\text{\rm \'et}}}

\newcommand{\Om}{{\Omega}}

\newcommand{\lci}{l.c.i.\ }

\newcommand{\uu}{\underline}
\newcommand{\oo}{\otimes}

\title{\bf  Algebraic Cobordism\vskip 6mm}

\author{M. Levine\vspace*{-0.5cm}\thanks{Department of Mathematics,
Northeastern University, Boston, MA 02115, USA. E-mail: marc@neu.edu}}

\date{\vspace{-8mm}}

\begin{document}
\maketitle

\thispagestyle{first} \setcounter{page}{57}

\begin{abstract}\vskip 3mm
Together with F. Morel, we have constructed in \cite{CR, Cobord1,
Cobord2} a theory of {\em algebraic cobordism}, an algebro-geometric version of the
topological theory of complex cobordism. In this paper, we give a survey of the
construction and main results of this theory; in the final section, we propose a
candidate for a theory of higher algebraic cobordism, which hopefully agrees with
the cohomology theory represented by the $\P^1$-spectrum $MGL$ in the Morel-Voevodsky
stable homotopy category.
\vskip 4.5mm


\end{abstract}

\vskip 12mm

\section{Oriented cohomology theories}\addtocounter{hack}{1}

\vskip-5mm \hspace{5mm}

Fix a field $k$ and let $\Sch_k$ denote the category of separated finite-type $k$-schemes. We let $\Sm_k$ be the
full subcategory of smooth quasi-projective $k$-schemes.

We have described in \cite{Cobord1} the notion of an {\em oriented cohomology theory}
on $\Sm_k$. Roughly speaking, such a theory $A^*$ consists of a contravariant functor
from $\Sm_k$ to graded  rings (commutative), which is also covariantly functorial for
projective equi-dimensional morphisms $f:Y\to X$ (with a shift in the grading):
\[
f_*:\A^*(Y)\to A^{*-\dim_XY}(X).
\]
The pull-back $g^*$ and push-forward $f_*$ satisfy a projection formula and
commute in transverse cartesian squares. If $L\to X$ is a line bundle with
zero-section $s:X\to L$, we have the {\em first Chern class} of $L$, defined by
\[
c_1(L):=s^*(s_*(1_X))\in A^1(X),
\]
where $1_X\in A^0(X)$ is the unit. $A^*$ satisfies the {\em projective bundle
formula}:
\begin{enumerate}
\item[(PB)] Let $\sE$ be a rank $r+1$ locally free coherent
sheaf on $X$, with projective bundle $q:\P(\sE)\to X$ and tautological quotient
invertible sheaf $q^*\sE\to \sO(1)$. Let $\xi=c_1(\sO(1))$. Then $A^*(\P(\sE))$ is a
free $A^*(X)$-module with basis $1,\xi,\ldots, \xi^r$.
\end{enumerate}
Finally, $A^*$ satisfies a homotopy property: if $p:V\to X$ is an affine-space bundle
(i.e., a torsor for a vector bundle over $X$), then $p^*:A^*(X)\to A^*(V)$ is an
isomorphism.

\begin{exs} (1) The theories $\CH^*$ and $H^{2*}_\et(-,\Z/n(*))$ on $\Sm_k$ (also with
$\Z_l(*)$ or $\Q_l(*)$ coefficients).

\smallskip
\noindent
(2) The theory $K_0[\beta,\beta^{-1}]$ on $\Sm_k$. Here $\beta$
is an indeterminant of degree $-1$, used to keep track of the relative dimension
when taking projective push-forward.
\end{exs}

\begin{rems}\label{rem:BMHom} (1) In \cite{Cobord2}, we consider a more general
(dual) notion, that of an {\em oriented Borel-Moore homology theory} $A_*$. Roughly,
this is a functor from a full subcategory of $\Sch_k$ to graded abelian groups,
covariant for projective maps, and contravariant (with a shift in the grading) for
local complete intersection morphisms. In addition, one has external products, and a
degree -1 Chern class endomorphism $\cn(L):A_*(X)\to A_{*-1}(X)$ for each line bundle
$L$ on
$X$, defined by $\cn(L)(\eta)=s^*(s_*(\eta))$, $s:X\to L$ the zero-section. As for an
oriented cohomology theory, there are various compatibilities of push-forward and
pull-back, and $A_*$ satisfies a projective bundle formula and a homotopy property.

This allows for a more general category of definition for $A_*$, e.g., the category
$\Sch_k$. As we shall see, the setting of Borel-Moore homology is often more natural
than cohomology. On $\Sm_k$, the two notions are equivalent: to pass from Borel-Moore
homology to cohomology, one re-grades by setting
$A^n(X):= A_{n-\dim_kX}(X)$ and uses the \lci pull-back for $A_*$ to give the
contravariant functoriality of $A^*$, noting that every morphism of smooth
$k$-schemes is an \lci morphism. We will state most of our results for
cohomology theories on $\Sm_k$, but they extend to the setting of Borel-Moore
homology on $\Sch_k$ (see \cite{Cobord2} for details).

\smallskip
\noindent
(2) Our notion of oriented cohomology is related to that of Panin \cite{Panin}, but
is not the same.
\end{rems}

\section{The formal group law}\addtocounter{hack}{1}

\vskip-5mm \hspace{5mm}

Let $A_*$ be an oriented
cohomology theory on $\Sm_k$. As noticed by Quillen \cite{Quillen}, a double
application of the projective bundle formula (PB) yields the isomorphism of rings
\[
A^*(k)[[u,v]]\cong \lim_{\substack{\leftarrow\\n,m}}A^*(\P^n\times\P^m),
\]
the isomorphism sending $u$ to $c_1(p_1^*\sO(1))$ and $v$ to $c_1(p_2^*\sO(1))$. The class of \linebreak
$c_1(p_1^*\sO(1)\otimes p_2^*\sO(1))$ thus gives a power series $F_A(u,v)\in A^*(k)[[u,v]]$ with
\[
c_1(p_1^*\sO(1)\otimes p_2^*\sO(1))=F_A(c_1(p_1^*\sO(1)),c_1(p_2^*\sO(1))).
\]
By the naturality of $c_1$, we have the identity for $X\in\Sm_k$ with line bundles
$L$,
$M$,
\[
c_1(L\otimes M)=F_A(c_1(L),c_1(M)).
\]
In addition, $F_A(u,v)=u+v\mod uv$, $F_A(u,v)=F_A(v,u)$, and
$F_A(F_A(u,v),w)=F_A(u,F_A(v,w))$. Thus, $F_A$ gives a formal group law
with coefficients in $A^*(k)$.

 \begin{rem} Note that $c_1:\Pic(X)\to A^1(X)$ is a
group homomorphism if and only if $F_A(u,v)=u+v$. If this is the case,
we call $A^*$ {\em ordinary}, if not, $A^*$ is {\em extraordinary}. If
$F_A(u,v)=u+v-\alpha uv$ with $\alpha$ a unit in $A^*(k)$, we call $A^*$ {\em
multiplicative and periodic}.
\end{rem}

\begin{exs} For $A^*=\CH^*$ or $H^{2*}$, $F_A=u+v$, giving examples of ordinary
theories. For the theory
$A=K_0[\beta,\beta^{-1}]$, $c_1(L)=(1-L^\vee)\beta^{-1}$, and $F_A(u,v)=u+v-\beta uv$,
giving an example of a multiplicative and periodic theory.
\end{exs}

\begin{rem} Let $\tilde{\L}^*=\Z[a_{ij}\ |\ i,j\ge 1]$, where we give
$a_{ij}$ degree $-i-j+1$, and let
$F\in\tilde{\L}^*[[u,v]]$ be the power series $F=u+v+\sum_{ij}a_{ij}u^iv^j$. Let
\[
\L^*=\tilde{\L}^*/F(u,v)=F(v,u), F(F(u,v),w)=F(u,F(v,w)),
\]
and let $F_\L\in\L^*[[u,v]]$ be the image of $F$.
Then $(F_\L, \L^*)$ is the universal commutative dimension 1 formal group; $\L^*$ is
called the {\em Lazard ring} ({\it cf.} \cite{Lazard}).

Thus, if $A^*$ is an oriented cohomology theory on $\Sm_k$, there is a canonical
graded ring homomorphism $\phi_A:\L^*\to A^*(k)$ with $\phi_A(F_\L)=F_A$.
\end{rem}

\section{Algebraic cobordism}\addtocounter{hack}{1}

\vskip-5mm \hspace{5mm}

The main result of \cite{Cobord1,
Cobord2} is

\begin{thm}\label{MainThm} Let $k$ be a field of characteristic zero.
\begin{enumerate}
\item There is a universal oriented Borel-Moore homology theory $\Omega_*$ on
$\Sch_k$. The restriction of $\Om_*$ to $\Sm_k$ yields the universal oriented
cohomology theory $\Om^*$ on $\Sm_k$.
\item The homomorphism $\phi_\Omega:\L^*\to\Omega^*(k)$ is an isomorphism.
\item Let $i:Z\to X$ be a closed imbedding with open complement $j:U\to X$. Then the
sequence
\[
\Omega_*(Z)\xrightarrow{i_*}\Omega_*(X)\xrightarrow{j^*}\Omega_*(U)\to 0
\]
is exact.
\end{enumerate}
\end{thm}

\noindent{\bf Idea of construction:}  We construct $\Omega_*(X)$ in steps; the construction is inspired by
Quillen's approach to complex cobordism \cite{Quillen}.
\begin{enumerate}
\item Start with {\em cobordism cycles} $(f:Y\to X, L_1,\ldots, L_r)$, with $Y\in
\Sm_k$ irreducible, $f:Y\to X$ projective and $L_1,\ldots, L_r$ line bundles on $Y$
(we allow $r=0$). We identify two cobordism cycles if there is an isomorphism
$\phi:Y\to Y'$, a permutation $\sigma$ and isomorphisms $L_j\cong
\phi^*L_{\sigma(j)}'$. Let
$\sZ_*(X)$ be the free abelian group on the cobordism cycles, graded by giving
$(f:Y\to X, L_1,\ldots, L_r)$ degree $\dim_kY-r$.
\item Let $\sR^{dim}(X)$ be the subgroup of $\sZ_*(X)$ generated by cobordism cycles
of the form $(f:Y\to X, \pi^*L_1,\ldots, \pi^*L_r, M_1,\ldots, M_s)$, where $\pi:Y\to
Z$ is a smooth morphism in $\Sm_k$, the $L_i$ are line bundles on $Z$, and
$r>\dim_kZ$. Let $\uu{\sZ}_*(X)=\sZ_*(X)/\sR^{dim}(X)$.
\item Add the Gysin isomorphism: If $L\to Y$ is a line bundle and $s:Y\to L$ is a
section transverse to the zero-section with divisor $i:D\to Y$, identify $(f:Y\to X,
L_1,\ldots, L_r, L)$ with $(f\circ i:D\to X, i^*L_1,\ldots, i^*L_r)$. We let
$\uu{\Om}_*(X)$ denote the resulting quotient of $\uu{\sZ}_*(X)$. Note that on
$\uu{\Om}_*(X)$ we have, for each line bundle $L\to X$, the {\em Chern class operator}
\begin{align*}
&\cn(L):\uu{\Om}_*(X)\to \uu{\Om}_{*-1}(X)\\
&(f:Y\to X, L_1,\ldots, L_r)\mapsto(f:Y\to X, L_1,\ldots, L_r,f^*L)
\end{align*}
as well as push-forward maps $f_*:\uu{\Om}_*(X)\to \uu{\Om}_*(X')$ for $f:X\to X'$
projective.
\item Impose the formal group law: Regrade $\L$ by setting
$\L_n:=\L^{-n}$. Let
$\Om_*(X)$ be the quotient of
$\L_*\oo\uu{\Om}_*(X)$ by the imposing the identity of maps $\L_*\oo\uu{\Om}_*(Y)\to
\L_*\oo\uu{\Om}_*(X)$
\[
(\id\oo f_*)\circ F_\L(\cn(L),\cn(M))=\id\oo (f_*\circ \cn(L\oo M))
\]
for $f:Y\to X$ projective, and $L, M$ line bundles on $Y$. Note that, having imposed
the relations in $\sR^{dim}$, the operators $\cn(L)$, $\cn(M)$ are locally nilpotent,
so the infinite series $F_\L(\cn(L),\cn(M))$ makes sense.
\end{enumerate}

As the notation suggests, the most natural construction of $\Om$ is as an oriented Borel-Moore homology theory
rather than an oriented cohomology theory; the tranlation to an oriented cohomology theory on $\Sm_k$ is given as
in remark~\ref{rem:BMHom}(1). The proof of theorem~\ref{MainThm} uses resolution of singularities \cite{Hironaka}
and the weak factorization theorem \cite{Abram} in an essential way.

\begin{rem}\label{rem:generators} In addition to the properties of $\Om_*$ listed in
theorem~\ref{MainThm},
$\Om_*(X)$ is generated by the classes of ``elementary" cobordism cycles $(f:Y\to X)$.
\end{rem}

\section{Degree formulas}\addtocounter{hack}{1}

\vskip-5mm \hspace{5mm}

In the paper \cite{Rost}, Rost made a
number of conjectures based on the theory of algebraic cobordism in the
Morel-Voevodsky stable homotopy category. Many of Rost's conjectures have been proved
by homotopy-theoretic means (see \cite{Borghesi}); our construction of algebraic
cobordism gives an alternate proof of these results, and settles many of the
remaining open questions as well. We give a sampling of some of these results.

\subsection{The generalized degree formula}

\vskip-5mm \hspace{5mm}

All the degree formulas follow from the  ``generalized degree formula". We first define the degree map
$\Om^*(X)\to \Om^*(k)$.

\begin{Def} Let $k$ be a field of characteristic zero and let $X$ be an irreducible
finite type $k$-scheme with generic point $i:x\to X$. For an element $\eta$ of
$\Om^*(X)$, define $\deg\eta\in\Om^*(k)$ to be the element mapping to $i^*\eta$ in
$\Om^*(k(x))$ under the isomorphisms $\Om^*(k)\cong \L^*\cong \Om^*(k(x))$ given by
theorem~\ref{MainThm}(2).
\end{Def}

\begin{thm}[(generalized degree formula)]\label{DecompThm} Let $k$ be a field of
characteristic zero. Let
$X$ be an irreducible finite type
$k$-scheme, and let
$\eta$ be in
$\Omega_*(X)$. Let $f_0:B_0\to X$ be a resolution of singularities of $X$, with $B_0$
quasi-projective over $k$. Then there are $a_i\in\Omega_*(k)$, and projective
morphisms $f_i:B_i\to X$ such that
\begin{enumerate}
\item  Each $B_i$ is in $\Sm_k$, $f_i:B_i\to f(B_i)$ is
birational and $f(B_i)$ is a proper closed subset of $X$ (for $i>0$).
\item $\eta-(\deg \eta)[f_0:B_0\to X]=\sum_{i=1}^r a_i[f_i:B_i\to X]$ in
$\Omega_*(X)$.
\end{enumerate}
\end{thm}

\begin{proof} It follows from the definitions of
$\Omega^*$ that we have
\[
\Omega^*(k(x))= \lim_{\substack{\to\\U}}\Omega^*(U),
\]
where the limit is over smooth dense open subschemes $U$ of $X$, and $\Omega^*(k(x))$
is the value at $\Spec k(x)$ of the functor $\Omega^*$ on finite type $k(x)$-schemes.
Thus, there is a smooth open subscheme $j:U\to X$ of $X$ such that
$j^*\eta=(\deg \eta)[\id_U]$  in $\Omega^*(U)$. Since $U\times_XB_0\cong U$, it
follows that $j^*(\eta-(\deg \eta)[f_0])=0$ in $\Omega^*(U)$.

Let $W=X\setminus U$. From the localization sequence
\[
\Omega_*(W)\xrightarrow{i_*}\Omega_*(X)\xrightarrow{j^*}\Omega_*(U)\to 0,
\]
we find an element $\eta_1\in\Omega_*(W)$ with $i_*(\eta_1)=\eta-(\deg \eta)[f_0]$,
and noetherian induction completes the proof.
\end{proof}

\begin{rem}\label{rem:DegForm} Applying theorem~\ref{DecompThm} to the class of a
projective morphism
$f:Y\to X$, with $X, Y\in \Sm_k$, we have the formula
\[
[f:Y\to X]-(\deg f)[\id_X]=\sum_{i=1}^ra_i[f_i:B_i\to X]
\]
in $\Om^*(X)$. Also, if $\dim_kX=\dim_kY$, $\deg f$ is the usual degree, i.e., the
field extension degree $[k(Y):k(X)]$ if $f$ is dominant, or zero if $f$ is not.
\end{rem}

\subsection{Complex cobordism}

\vskip-5mm \hspace{5mm}

For a differentiable manifold $M$, one has the complex cobordism ring $MU^*(M)$. Given an embedding $\sigma:k\to
\C$ and an $X\in \Sm_k$, we let $X^\sigma(\C)$ denote the complex manifold associated to the smooth $\C$-scheme
$X\times_k\C$. Sending $X$ to $MU^{2*}(X^\sigma(\C))$ defines an oriented cohomology theory on $\Sm_k$; by the
universality of $\Om^*$, we have a natural homomorphism
\[
\Re_\sigma:\Omega^*(X)\to MU^{2*}(X^\sigma(\C)).
\]

Now, if $P=P(c_1,\ldots, c_d)$ is a degree $d$ (weighted) homogeneous
polynomial, it is known that the operation of sending a smooth compact
$d$-dimensional complex manifold $M$ to the Chern number $\deg(P(c_1,\ldots,
c_d)(\Theta_M))$ (where $\Theta_M$ is the complex tangent bundle) descends to a
homomorphism $MU^{-2d}\to\Z$. Composing with $\Re_\sigma$, we have the homomorphism
$P:\Om^{-d}(k)\to \Z$. If $X$ is smooth and projective of dimension $d$ over $k$, we
have $P([X])=\deg(P(c_1,\ldots,c_d)(\Theta_{X^\sigma(\C)}))$;
$P([X])$ is in fact independent of the choice of embedding $\sigma$.

Let $s_d(c_1,\ldots, c_d)$ be the polynomial
which corresponds to $\sum_i\xi_i^d$, where $\xi_1,\ldots$ are the Chern roots. The
following divisibility is known (see \cite{Adams}):  if $d=p^n-1$ for some prime $p$,
and $\dim X=d$, then $s_d(X)$ is divisible by $p$.

In addition, for integers $d=p^n-1$ and  $r\ge 1$, there are mod $p$ characteristic
classes $t_{d,r}$, with $t_{d,1}=s_d/p\mod p$. The
$s_d$ and the
$t_{d,r}$ have the following properties:
\begin{equation}\label{Properties}
\end{equation}
\begin{enumerate}
\item $s_d(X)\in p\Z$ is defined for $X$ smooth and projective of dimension
$d=p^n-1$. $t_{d,r}(X)\in \Z/p$ is defined for $X$ smooth and projective of dimension
$rd=r(p^n-1)$.
\item $s_d$ and $t_{d,r}$ extend to homomorphisms
$s_d:\Omega^{-d}(k)\to p\Z$, $t_{d,r}:\Omega^{-rd}(k)\to \Z/p$.
\item If $X$ and $Y$ are smooth projective varieties with $\dim X, \dim Y>0$, $\dim
X+\dim Y=d$, then
$s_d(X\times Y)=0$.
\item If $X_1,\ldots, X_s$ are smooth projective varieties with $\sum_i\dim X_i=rd$,
then $t_{d,r}(\prod_iX_i)=0$ unless $d|\dim X_i$ for each $i$.
\end{enumerate}

We can now state Rost's degree formula and the higher degree formula:

\begin{thm}[(Rost's degree formula)]\label{DegFormula} Let $f:Y\to X$ be a morphism of
smooth projective
$k$-schemes of dimension $d$, $d=p^n-1$ for some prime $p$. Then there is a
zero-cycle $\eta$ on
$X$ such that
\[
s_d(Y)-(\deg f)s_d(X)=p\cdot \deg(\eta).
\]
\end{thm}

\begin{thm}[(Rost's higher degree formula)]\label{GenDegFormula} Let $f:Y\to X$ be a
morphism of smooth projective
$k$-schemes of dimension $rd$, $d=p^n-1$ for some prime $p$. Suppose that $X$ admits
a sequence of surjective morphisms
\[
X=X_0\to X_1\to\ldots\to X_{r-1}\to X_r=\Spec k,
\]
such that:
\begin{enumerate}
\item $\dim X_i=d(r-i)$.
\item Let $\eta$ be a zero-cycle on
$X_i\times_{X_{i+1}}\Spec k(X_{i+1})$. Then $p|\deg(\eta)$.
\end{enumerate}
Then
\[
t_{d,r}(Y)=\deg(f)t_{d,r}(X).
\]
\end{thm}

\begin{proof} These two theorems follow easily from the generalized degree formula.
Indeed, for theorem~\ref{DegFormula}, take the identity of remark~\ref{rem:DegForm}
and push forward to $\Om^*(k)$. Using remark~\ref{rem:generators}, this gives the
identity
\[
[Y]-(\deg f)[X]=\sum_{i=1}^r m_i[A_i\times B_i]
\]
in $\Omega^*(X)$, for smooth, projective $k$-schemes $A_j$, $B_j$, and integers
$m_j$, where each $B_i$ admits a projective morphism $f_i:B_i\to X$ which is
birational to its image and not dominant. Since $s_d$ vanishes on non-trivial
products, the only relevant part of the sum involves those $B_j$ of dimension zero;
such a $B_j$ is identified with the closed point $b_j:=f_j(B_j)$ of $X$. Applying
$s_d$, we have
\[
s_d(Y)-\deg(f) s_d(X)=\sum_jm_js_d(A_j)\deg_k(b_j).
\]
Since $s_d(A_j)=pn_j$ for suitable integers $n_j$, we
have
\[
s_d(Y)-\deg(f)s_d(X)=p\deg(\sum_j m_jn_jb_j).
\]
Taking $\eta=\sum_j m_jn_jb_j$ proves theorem~\ref{DegFormula}.

The proof of theorem~\ref{GenDegFormula} is similar: Start with the
decomposition of $[f:Y\to X]-(\deg f)[\id_X]$ given by remark~\ref{rem:DegForm}. One
then decomposes the maps $B_i\to X=X_0$ further by pushing forward to $X_1$ and using
theorem~\ref{DecompThm}. Iterating down the tower gives the identity in
$\Omega_*(k)$
\[
[Y]-(\deg f)[X]=\sum_i m_i[B^i_0\times \ldots\times B^i_r];
\]
the condition (2) implies that, if $d|\dim_kB^i_j$ for all $j=0,\ldots, r$, then
$p|m_j$. Applying $t_{d,r}$ and using the property \eqref{Properties}(4) yields the
formula.
\end{proof}

\section{Comparison results}\addtocounter{hack}{1}

\vskip-5mm \hspace{5mm}

Suppose we have a formal group $(f,R)$, giving the canonical homomorphism
$\phi_f:\L^*\to R$. Let $\Omega^*_{(f,R)}$ be the functor
\[
\Omega^*_{(f,R)}(X)=\Omega^*(X)\otimes_{\L^*}R,
\]
where $\Omega^*(X)$ is an $\L^*$-algebra via the homomorphism
$\phi_\Omega:\L^*\to\Omega^*(k)$.  The universal property of $\Omega^*$ gives
the analogous universal property for $\Omega^*_{(f,R)}$.

In particular, let $\Omega^*_+$ be the theory with $(f(u,v),R)=(u+v,\Z)$, and let
$\Omega^*_\times$ be the theory with $(f(u,v),R)=(u+v-
\beta uv,\Z[\beta,\beta^{-1}])$. We thus have the canonical natural transformations of
oriented theories on $\Sm_k$
\begin{equation}\label{Nat}
\Omega^*_+\to \CH^*;\quad \Omega^*_\times\to K_0[\beta,\beta^{-1}].
\end{equation}

\begin{thm}\label{} Let $k$ be a field of characteristic zero. The natural
transformations
\eqref{Nat} are isomorphisms, i.e., $\CH^*$ is the universal ordinary oriented
cohomology theory and $K_0[\beta,\beta^{-1}]$ is the universal multiplicative
and periodic theory.
\end{thm}

\begin{proof} For $\CH^*$, this uses localization, theorem~\ref{DecompThm} and
resolution of singularities. For $K_0$, one writes down an integral Chern character,
which gives the inverse isomorphism by the Grothendieck-Riemann-Roch theorem.
\end{proof}

\section{Higher algebraic cobordism}

\vskip-5mm \hspace{5mm}

The cohomology theory represented by the
$\P^1$-spectrum
$MGL$ in the Morel-Voevod\-sky
$\A^1$-stable homotopy category \cite{MorelVoev, Voev} gives perhaps the most natural
algebraic analogue of complex cobordism. By universality,
$\Om^n(X)$ maps to
$MGL^{2n,n}(X)$; to show that this map is an isomorphism, one would like to give a map
in the other direction. For this, the most direct method would be to extend $\Om^*$ to
a theory of higher algebraic cobordism; we give one possible approach to this
construction here.

The idea is to repeat the construction of $\Om_*$, replacing abelian
groups with symmetric monoidal categories throughout. Comparing with the
$Q$-construction, one sees that the cobordism cycles in $\sR^{dim}(X)$ should be
homotopic to zero, but not canonically so. Thus, we cannot impose this relation
directly, forcing us to modify the group law by taking a limit.

Start with the category $\widetilde{\sZ}(X)_0$, with objects $(f:Y\to X, L_1,\ldots,
L_r)$, where
$Y$ is irreducible in $\Sm_k$, $f$ is projective, and the $L_i$ are line bundles on
$Y$.  A morphism $(f:Y\to X, L_1,\ldots, L_r)\to (f':Y'\to X, L_1',\ldots, L_r')$
in $\widetilde{\sZ}(X)_0$  consist of a tuple $(\phi, \psi_1,\ldots, \psi_r,\sigma)$,
with $\phi:Y\to Y'$ an isomorphism over $X$, $\sigma$ a permutation, and
$\psi_j:L_j\to \phi^*L_{\sigma(j)}'$ an isomorphism of line bundles on $Y$. Form the
category $\widetilde{\sZ}(X)$ as the symmetric monoidal category freely generated by
$\widetilde{\sZ}(X)_0$; grade $\widetilde{\sZ}(X)$ by letting $\widetilde{\sZ}_n(X)$
be the full symmetric monoidal subcategory generated by the $(f:Y\to X, L_1,\ldots,
L_r)$ with $n=\dim_kY-r$.

Next, form $\widetilde{\uu{\Om}}(X)$ by adjoining (as a symmetric monoidal
category) an isomorphism
$\gamma_{L,s}:(f\circ i:D\to X, i^*L_1,\ldots, i^*L_r)\to (f:Y\to X, L_1,\ldots, L_r,
L)$ for each section $s:Y\to L$ transverse to the zero-section with divisor $i:D\to
X$. Given a morphism
$\tilde{\phi}:=(\phi,\ldots):(f:Y\to X, L_1,\ldots, L_r, L)\to
(f':Y'\to X, L_1',\ldots, L_r', L')$ (with $L\cong \phi^*L'$ via $\tilde\phi$), let
$i':D'\to Y'$ be the map induced by
$\phi$,
$s':Y'\to L'$ the section induced by $s$, and
\[
\psi^D:(f\circ i:D\to X,i^*L_1,\ldots,i^*L_r)\to (f'\circ i':D'\to X,i^{\prime*}L_1',\ldots,i^{\prime*}L_r'),
\]
the morphism induced by $\psi$. We impose the relation
$\psi\circ\gamma_{L,s}=\gamma_{L',s'}\circ\psi^D$. Finally, for line bundles $L, M$
with smooth transverse divisors $i_D:D\to Y$,
$i_E:E\to Y$ defined by sections
$s:Y\to L$, $t:Y\to M$, respectively, we impose the relation
$\gamma_{L,s}\circ\gamma_{i_D^*M,i_D^*t}=
\gamma_{M,t}\circ\gamma_{i_E^*L,i_E^*s}$.
The grading on
$\widetilde{{\sZ}}(X)$ extends to one on
$\widetilde{\uu{\Om}}(X)$.

Given $g:X\to X'$ projective, we have the functor
$g_*:\widetilde{\uu{\Om}}(X)\to
\widetilde{\uu{\Om}}(X')$, similarly, given a smooth morphism $h:X\to X'$, we have
the functor $h^*:\widetilde{\uu{\Om}}(X')\to \widetilde{\uu{\Om}}(X)$. Given a line
bundle $L$ on $X$, we have the natural transformation $\cn(L)$ sending $(f:Y\to
X,L_1,\ldots, L_r)$ to $(f:Y\to
X,L_1,\ldots, L_r, f^*L)$.

Now let $\sC$ be a symmetric monoidal category such that all morphisms are
isomorphisms, and let
$R$ be a ring, free as a $\Z$-module. One can define a symmetric monoidal category
$R\otimes_\N \sC$ with a symmetric monoical functor $\sC\to R\otimes_\N \sC$ which is
universal for symmetric monoidal functors $\sC\to \sC'$ such that $\sC'$ admits an
action of $R$ via natural transformations. In case $R=\Z$, $\Z\otimes_\N\sC$ is the
standard group completion $\sC^{-1}\sC$. In general, if $\{e_\alpha\ |\ \alpha\in
A\}$ is a $\Z$-basis for $R$, then
\[
R\otimes_\N \sC=\coprod_\alpha \sC^{-1}\sC,
\]
with the $R$-action given by expressing $\times x:R\to R$ in terms of the basis
$\{e_\alpha\}$.

For each integer $n\ge0$, let $\L^{(n)}_*$ be the quotient of $\L_*$ by the ideal of
elements of degree $>n$. We thus have the formal group $(F_{\L^{(n)}},\L^{(n)}_*)$.

 We form the
category $\L^{(n)}\otimes_\N\widetilde{\uu{\Om}}(X)$, which we grade by total
degree. For each
$f:Y\to X$ projective, with
$Y\in\Sm_k$, and line bundles $L$, $M$, $L_1,\ldots, L_r$ on $Y$, we adjoin an
isomophism $\rho_{L,M}$
\[
f_*(F_{\L^{(n)}}(\cn(L),\cn(M))(\id_Y,L_1,\ldots, L_r))\xrightarrow{\sim}
f_*(\id\otimes\cn(L\otimes M)(\id_Y,L_1,\ldots, L_r)).
\]
We impose the condition of naturality with respect to the maps in
$\L^{(n)}\otimes_\N\widetilde{\uu{\Om}}_n(Y)$, in the evident sense; the Chern class
transformations extend in the obvious manner.

We impose the following commutativity condition: We have the  evident
isomorphism $t_{L,M}:F_{\L^{(n)}}(\cn(L),
\cn(M))\to  F_{\L^{(n)}}(\cn(M),\cn(L))$  of natural
transformations, as well as $\tau_{L,M}:\cn(L\oo M)\to \cn(M\oo L)$,
the isomorphism induced by the symmetry $L\oo M\cong M\oo L$. Then we impose
the identity $\tau_{L,M}\circ \rho_{L,M}=\rho_{M,L}\circ t_{L,M}$.
We impose a similar identity between the associativity of the formal group law and
the associativity of the tensor product of line bundles.

We also adjoin $a\cdot\tau_{L,M}$ for all $a\in\L^{(n)}$, with similar
compatibilities as above,  respecting the $\L^{(n)}$-action and sum. This forms
the symmetric monoidal category
$\widetilde{\Om}^{(n)}(X)$, which inherits a grading from $\widetilde{\uu{\Om}}(X)$.
We have the inverse system of graded symmetric monoidal categories:
\[
\ldots\to\widetilde{\Om}^{(n+1)}(X)\to \widetilde{\Om}^{(n)}(X)\to\ldots\ .
\]

\begin{Def} Set $\Om_{m,r}^{(n)}(X):=\pi_r(B\widetilde{\Om}_m^{(n)}(X))$ and
$\Om_{m,r}(X):=\displaystyle\lim_{\substack{\leftarrow\\n}}\Om_{m,r}^{(n)}(X)$.
\end{Def}

At present, we can only verify the following:

\begin{thm} There is a natural isomorphism $\Om_{m,0}(X)\cong \Om_m(X)$.
\end{thm}

\begin{proof} First note that $\pi_0(\widetilde{\sZ}_m(X))$ is a
commutative monoid with group completion $\sZ_m(X)$. Next, the natural map
$\pi_0(\widetilde{\uu{\Om}}_*(X))^+\to \uu{\Om}_*(X)$
is surjective with kernel generated by the classes generating $\sR^{dim}(X)$. Given
such an element $\psi:=(f:Y\to X, \pi^*L_1,\ldots, \pi^*L_r, M_1,\ldots, M_s)$, with
$\pi:Y\to Z$ smooth, and $r>\dim_kZ$, suppose that the $L_i$ are very ample. We may
then choose sections $s_i:Z\to L_i$ with divisors $D_i$ all intersecting
transversely. Iterating the isomorphisms $\gamma_{L_i, s_i}$ gives a path from
$\psi$ to 0 in $B\widetilde{\uu{\Om}}_*(X)$. Passing to
$B\widetilde{\Om}_m^{(n)}(X)$, the group law allows us to replace an arbitrary line
bundle witha difference of very ample ones, so all the classes of this form go to
zero in $\Om_{m,0}^{(n)}(X)$. This shows that the natural map
\[
\Om_{m,0}^{(n)}(X)\to (\L^{(n)}\oo_\L\Om_*(X))_m
\]
is an isomorphism. Since $(\L^{(n)}\oo_\L\Om_*(X))_m=\Om_m(X)$ for $m\ge n$, we are
done.
\end{proof}

The categories
$\widetilde{\Om}^{(n)}_m(X)$ are covariantly functorial for projective maps,
contravariant for smooth maps (with a shift in the grading) and have first
Chern class natural transformations $\cn(L):\widetilde{\Om}^{(n)}_m(X)\to
\widetilde{\Om}^{(n)}_{m-1}(X)$ for $L\to X$ a line bundle.

We conjecture that the inverse system used to define $\Om_{m,r}(X)$ is eventually
constant for all $r$, not just for $r=0$. If this is true, it is reasonable to define
the space $B\widetilde{\Om}_m(X)$ as the homotopy limit
\[
B\widetilde{\Om}_m(X):=\holim_n B\widetilde{\Om}_m^{(n)}(X).
\]
One would then have $\Om_{m,r}(X)=\pi_r(B\widetilde{\Om}_m(X),0)$ for all $m,r$;
hopefully the properties of $\Om_*$ listed in theorem~\ref{MainThm} would then
generalize into properties of the spaces $B\widetilde{\Om}_m(X)$.

\label{lastpage}

\end{document}